\documentclass[twocolumn]{autart}
\usepackage{amsfonts,amsmath,amssymb}
\usepackage{lineno,hyperref}
\usepackage{subfigure}
\modulolinenumbers[5]
\usepackage{graphicx}
\graphicspath{{./figures_eps/}}
\usepackage{cases}
\usepackage{comment}
\newtheorem{mydef}{Definition}
\newtheorem{mythm}{Theorem}
\newtheorem{mylem}{Lemma}
\newtheorem{myas}{Assumption}
\newtheorem{myrem}{Remark}

\newtheorem{myprob}{Problem}

\newcommand{\cal}{\mathcal}

\newcommand{\rfig}[1]{Fig.\,\ref{#1}} 
 
\newcommand{\req}[1]{\eqref{#1}}

\newcommand{\rlem}[1]{Lemma\,\ref{#1}}
\newcommand{\rdef}[1]{Definition\,\ref{#1}}
\newcommand{\rrem}[1]{Remark\,\ref{#1}}
\newcommand{\ras}[1]{Assumption\,\ref{#1}}
\newcommand{\rthm}[1]{Theorem\,\ref{#1}}
\newcommand{\rprob}[1]{Problem\,\ref{#1}}
\renewcommand{\labelenumi}{(\roman{enumi})}
\renewcommand\thefootnote{\arabic{footnote})}
\newcommand{\qedwhite}{\hfill \ensuremath{\Box}}

\begin{document}

\begin{frontmatter}

\title{Event-Triggered Intermittent Sampling for \\ Nonlinear Model Predictive Control\thanksref{footnoteinfo}}
\thanks[footnoteinfo]{This work was supported by the Swedish Research Council (VR), Knut och Alice Wallenberg foundation (KAW), and the H2020 ERC Starting Grant BUCOPHSYS.}
\author[Second]{Kazumune Hashimoto}\ead{kazumune.hashimoto@z5.keio.jp},
\author[Second]{Shuichi Adachi}\ead{adachi@appi.keio.ac.jp},
\author[Third]{Dimos V. Dimarogonas}\ead{dimos@kth.se}
\address[Second]{School of Applied Physics and Physico informatics, Keio University, 
   Yokohama, Japan}\address[Third]{School of Electrical Engineering, KTH Royal Institute of Technology, 10044 Stockholm, Sweden.}

\begin{abstract}
In this paper, we propose a new aperiodic formulation of model predictive control for nonlinear continuous-time systems. Unlike earlier approaches, we provide event-triggered conditions \textit{without} using the optimal cost as a Lyapunov function candidate. Instead, we evaluate the time interval when the optimal state trajectory enters a local set around the origin. The obtained event-triggered strategy is more suitable for practical applications than the earlier approaches in two directions. First, it does not include parameters (e.g., Lipschitz constant parameters of stage and terminal costs) which may be a potential source of conservativeness for the event-triggered conditions. Second, the event-triggered conditions are necessary to be checked only at certain sampling time instants, instead of continuously. This leads to the alleviation of the sensing cost and becomes more suitable for practical implementations under a digital platform. The proposed event-triggered scheme is also validated through numerical simulations.
\end{abstract}


\end{frontmatter}


\section{Introduction}
Event-Triggered Control (ETC) and Self-Triggered Control (STC) have been active areas of research in the community of Networked Control Systems (NCSs), due to their potential advantages over the typical time-triggered controllers \cite{heemels2012a}. 
In contrast to the time-triggered case where the control signals are executed periodically, ETC and STC trigger the executions based on the violation of certain prescribed control performances, see e.g., \cite{heemels2011a,lemmon2009a}. 

In another line of research, Model Predictive Control (MPC) has been one of the most popular control strategies applied in a wide variety of applications. MPC plays an important role when several constraints, such as actuator or physical limitations, need to be explicitly taken into account. 
The basic idea of MPC is to obtain the current control action by solving the Optimal Control Problem (OCP) online, based on the knowledge of current state measurement and future behavior prediction through the dynamics. 

The application of ETC and STC framework to MPC, generally known as Event-Triggered MPC (ETMPC) and Self-triggered MPC (STMPC), is particular of importance as it potentially alleviates a computational load by reducing the amount of solving OCPs. In ETMPC and STMPC, the OCPs are solved only when some events, generated based on certain control performance criteria, are triggered. These strategies have received an increased attention in recent years; most of the works focus on discrete-time systems, see e.g., \cite{schedule5,evmpc_linear10,evmpc_linear4,evmpc_linear5,evmpc_linear6,evmpc_nonlinear4,evmpc_nonlinear5,hashimoto2015c}, and some results include for the continuous-time case, see e.g., \cite{hashimoto2016a,evmpc_linear9,evmpc_linear3} for linear systems and  \cite{evmpc_nonlinear1,evmpc_nonlinear6,evmpc_nonlinear9,evmpc_nonlinear2,hashimoto2015a,hashimoto2017a} for nonlinear systems. 
In this paper, we are particularly interested in the case of nonlinear continuous-time systems. 
Among the afore-cited papers for nonlinear continuous-time systems, the results can be further divided into two categories, depending on whether disturbances are taken into account; see \cite{evmpc_nonlinear6} for the disturbance-free case and \cite{evmpc_nonlinear1,evmpc_nonlinear2,evmpc_nonlinear9,hashimoto2015a,hashimoto2017a} for the presence of disturbance case. 
In \cite{evmpc_nonlinear6}, an event-triggered MPC strategy has been proposed for nonlinear systems with no disturbances. While a delay compensation strategy has been developed to tackle uncertainties for networked control systems, an explicit form of the event-triggered condition is not provided and beyond the scope of that paper. In \cite{evmpc_nonlinear9}, a self-triggered strategy is proposed for general nonlinear systems with additive disturbances. The self-triggered condition was derived based on the optimal cost regarded as an ISS Lyapunov function candidate. In \cite{evmpc_nonlinear2}, an event-triggered strategy has been proposed for general nonlinear systems with additive bounded disturbances. When deriving the event-triggered strategy, an additional state constraint is imposed such that the optimal cost as a Lyapunov function candidate is decreasing. In \cite{hashimoto2017a}, a self-triggered strategy was provided for nonlinear input affine systems based on the optimal cost as a Lyapunov function candidate. In the approach, an additional way to discretize an optimal control trajectory into several control {samples} was provided so that these can be transmitted to the plant over the network channels. 

In this paper, we propose a new event-triggered formulation of MPC for nonlinear continuous-time systems with additive bounded disturbances. The main novelty of the proposed framework with respect to earlier results in this category (\cite{evmpc_nonlinear1,evmpc_nonlinear2,evmpc_nonlinear9,hashimoto2015a,hashimoto2017a}), is that the event-triggered condition is derived based on a new stability theorem, which does not evaluate the optimal cost as a Lyapunov function candidate. In the stability derivations, we instead evaluate the \textit{time interval}, when the optimal state trajectory enters a local region around the origin. By guaranteeing that this time interval becomes smaller as the OCP is solved, it is ensured that the state enters a prescribed set in finite time. 

The derivation of the new stability is motivated by the fact that the earlier event-triggered strategies may include Lipschitz constant parameters for the stage and terminal cost (see e.g., \cite{evmpc_nonlinear1,hashimoto2017a}). 
When standard quadratic costs are utilized, the corresponding Lipschitz parameters are characterized by the {maximum} distance of the state from the origin \cite{evmpc_nonlinear1}, and the triggering condition becomes largely affected by the state domain considered in the problem formulation. That is, as a larger state domain is considered, the event-triggered condition may become more conservative. 
Depending on the problem formulation, therefore, it may not be desirable to include these parameters in the event-triggered condition. Since our approach does not evaluate the optimal cost as a Lyapunov functon candidate,  the corresponding event-triggered conditions do not include such unsuitable parameters even though quadratic cost functions are used. We will also illustrate through a simulation example that the proposed approach attains much less conservative result than our previous result presented in \cite{hashimoto2017a}.

As another contribution of this paper with respect to the 
afore-cited papers of ETMPC for continuous-time systems (including the linear case), we will additionally incorporate {Periodic} Event-Triggered Control (PETC) framework \cite{heemels2013a}. In PETC, triggering conditions are evaluated only at certain sampling time instants, instead of continuously. This approach has certain advantages over the existing ETMPC strategies, since it alleviates a sensing load to evaluate the event-triggered conditions and becomes more suitable to be implemented under digital platforms. 
In the general PETC framework, the sampling time to evaluate the event-triggered condition is constant for all update times \cite{heemels2013a}. In our proposed approach, on the other hand, the sampling time is selected in an adaptive way; for each time of solving OCP, the controller adaptively determines the sampling time to check the event-triggered condition, such that the desired control performance can be guaranteed.

This paper is organized as follows. 
In Section 2, the optimal control problem is formulated. In Section 3, feasibility of the OCP is analyzed. In Section 4, our main proposed algorithm is presented, and the stability is shown in Section 5. A simulation example validates our proposed method in Section 6. We finally conclude in Section 7. 

\noindent
\textbf{Notations.} Let $\mathbb{R}$, $\mathbb{R}_{> 0}$, $\mathbb{R}_{\geq 0}$,  $\mathbb{N}_{\geq 0}$, $\mathbb{N}_{\geq 1}$ be the real, positive real, non-negative real, non-negative integers and positive integers, respectively. For a given matrix $Q$, we use $Q \succ 0$ to denote that the matrix $Q$ is positive definite. The notation $\lambda_{\min} (Q)$ is used to denote the minimal eigenvalue of the matrix $Q$. We denote $||x||$ as the Euclidean norm of vector $x$, and $||x||_P$ as a weighted norm of vector $x$, i.e., $||x||_P = \sqrt{x^\mathsf{T} P x}$. Given a compact set $\Phi \subseteq \mathbb{R}^n $, we denote by $\partial \Phi$ the boundary of $\Phi$. The function $f : \mathbb{R}^n \times  \mathbb{R}^m \rightarrow \mathbb{R}^n$ is called Lipschitz continuous in $\mathbb{R}^n$ with a weighted matrix $P$, if there exists $0\leq L_{f} < \infty$ such that $||f (x_1, u) - f (x_2, u) ||_P \leq L_{f}  ||x_1 - x_2||_P$, $\forall x_1, x_2 \in \mathbb{R}^n$, $\forall u \in \mathbb{R}^m$. 

\section{Problem formulation}
\subsection{Dynamics and optimal control problem}
In this section the problem formulation is defined. 
We consider to apply MPC to the following nonlinear systems with additive disturbances:
\begin{equation}\label{sys1} 
\dot{{x}}(t) = f({x}(t) , u(t)) + w(t), \ \ t \geq t_0,
\end{equation}
where ${x}(t)\in \mathbb{R}^n$ is the state, $u(t)\in \mathbb{R}^m$ is the control input, $w(t)\in \mathbb{R}^n$ is an additive bounded disturbance, and $t_0 \in \mathbb{R}$ denotes the initial time. The control input $u$ and the disturbance $w$ are assumed to satisfy the following constraints:
\begin{equation}\label{constraintuw}
u(t) \in {\mathcal U}\subseteq \mathbb{R}^m, \ w(t)\in {\mathcal W}\subseteq \mathbb{R}^n, \ \ \forall t \geq t_0.
\end{equation}
Regarding the constraint \req{constraintuw} and the plant model \req{sys1}, we make the following standard assumptions \cite{Chen1998a}:  
\begin{myas}\label{as1}
(i) The constraint sets ${\mathcal U}$ and ${\mathcal W}$ are compact, convex and $0 \in {\mathcal U}$; (ii) 
the function $f : \mathbb{R}^n \times  \mathbb{R}^m \rightarrow \mathbb{R}^n$ is twice continuously differentiable, and $f (0,0)=0$; (iii) the system \req{sys1} has a unique, absolutely continuous solution for any initial state $x(t_0)$ and any piecewise continuous control and disturbance $u: [t_0, \infty) \rightarrow {\cal U}$, $w: [t_0, \infty) \rightarrow {\cal W}$;
(iv) for the linearized system around the origin with no disturbances, i.e., 
$\dot{x}(t) = A_f x(t) + B_f u(t)$, 
where $A_f = \partial f /\partial x (0, 0)$ and $B_f = \partial f/ \partial u (0, 0)$, the pair $(A_f, B_f)$ is stabilizable. 
\end{myas}
Let $t_k$, $k\in \mathbb{N}_{\geq 0}$ be the update time instants when OCPs are solved, and let $\Delta_k = t_{k+1} - t_k$ be the inter-event times. At $t_k$, the controller solves an OCP based on the state measurement $x(t_k)$ and the predictive behavior of the systems described by \req{sys1}. In this paper, we consider the following cost to be minimized: 
\begin{equation}\label{cost}
\begin{array}{lll}
J(x(t_k), u(\cdot )) ={\displaystyle \int}^{t_k+T_k}_{t_k} || \hat{x}(\xi ) ||^2 _Q + || u(\xi )||^2 _R {\rm d}\xi,
\end{array}
\end{equation}
where $Q= Q^\mathsf{T} \succ 0$, $R = R^\mathsf{T} \succ 0$ and $T_k >0$ is the prediction horizon. $\hat{x}(\xi )$ denotes the nominal trajectory of \req{sys1} given by 
$\dot{\hat{x}} (\xi ) = f( \hat{x} ( \xi), u (\xi ) )$
for all $\xi \in [t_k , t_k + T_k]$ with $\hat{x} (t_k) = x(t_k)$. 
Here, the prediction horizon $T_k$ is not constant but is adaptively selected such that it is strictly decreasing. More characterization of $T_k$ is provided in this section when formulating the OCP. 

The following lemma states that there exists a stabilizing, state feedback controller in a prescribed local set around the origin:  
\begin{mylem}\label{lem1}
Suppose that \ras{as1} holds. 
Then, there exists a positive constant $0<\varepsilon <\infty$, a matrix $P_f =P ^\mathsf{T} _f \succ 0$, and a local controller $\kappa(x) = K x \in {\mathcal U}$, satisfying 
\begin{equation}\label{robust_terconst}
\frac{{\partial} V_f}{{\partial }x}\ f (x, \kappa(x) ) \leq - \cfrac{1}{2}\ {x}^\mathsf{T} (Q+K^\mathsf{T} R K) {x} 
\end{equation}
for all $x \in \Phi$, where $V_f (x) = {x}^\mathsf{T} P_f {x}$ and $\Phi = \{ x \in \mathbb{R}^{n} : V_f (x) \leq \varepsilon^2 \}$. Furthermore, $\Phi$ is a positive invariant set for the system \req{sys1} with $\kappa(x) = K x \in {\mathcal U}$, if the disturbance $w$ satisfies $||w||_{P_f} \leq \hat{w} _{\rm max}$ with $\hat{w} _{\rm max} = \varepsilon\lambda_{\min} (\hat{Q}_P ) /4$ and $\hat{Q}_{P} =P^{-1/2} _f ( Q+K^\mathsf{T} R K )P^{-1/2} _f $.
\end{mylem}
The proof is obtained by extending \textit{Lemma 1} in \cite{Chen1998a} and is given in the Appendix. 
\begin{mydef}[Control Objective of MPC]
The control objective of MPC is to steer the state $x$ to the local region $\Phi$ in finite time. 
\end{mydef}

In this paper, we consider that the control law switches from applying MPC to the utilization of the local controller $\kappa$, as soon as the state enters $\Phi$. 
This switching control law is referred to as \textit{dual mode MPC}, which is adopted in many works in the literature \cite{michalska1993}. Note that if the plant is controlled over a network \footnote{Note that our problem formulation is not limited to NCSs. Since the objective of this paper is to reduce the computation load of solving OCPs, applying our approach is still useful even for the case when the plant is not controlled over a network.}, applying the local controller $\kappa (x)$ may require a \textit{continuous} control update and may not be suitable under limited communication capabilities. One way to avoid this issue is to apply the local controller in a \textit{sample-and-hold fashion}, i.e., $u(t) = \kappa (x(t_k)),\ t \in [t_k, t_k + \delta]$. Here, $0< \delta < \infty$ can be chosen small enough such that asymptotic stability is still guaranteed, see \cite{magni2004a} for a detailed analysis. Also, please see \textit{Remark 2} in \cite{hashimoto2017a} for yet another way to avoid the problem of such continuous requirement. 

\begin{figure}[tbp]
  \begin{center}
   \includegraphics[width=7.0cm]{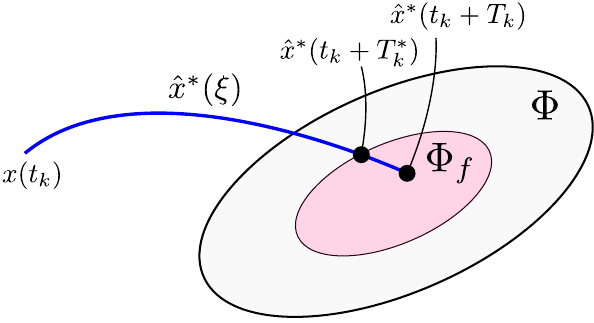}
   \caption{Graphical representation of the two regions $\Phi$ (grey region), $\Phi_f$ (red region) and the optimal state trajectory $\hat{x}^*$ (blue solid line). $T^* _k$ denotes the time interval to reach $\Phi_f$. } 
   \label{terregion_robust}
  \end{center}
 \end{figure}

Based on the local set $\Phi$, we further define the restricted set $\Phi_f$ given by ${\Phi}_f = \{ x \in \mathbb{R}^{n} : V_f (x) \leq \varepsilon^2 _f \}$,
where $0< \varepsilon_f <\varepsilon$. 
Since $\varepsilon_f <\varepsilon$, the set ${\Phi}_f $ is contained in $\Phi$, i.e., ${\Phi}_f \subset \Phi$. An example of these two regions is illustrated in \rfig{terregion_robust}. 
\begin{myas}\label{lipschitzf}
The nonlinear function $f : \mathbb{R}^n \times  \mathbb{R}^m \rightarrow \mathbb{R}^n$ is Lipschitz continuous in $\mathbb{R}^n$ with the weighted matrix $P_f$, with the Lipschitz constant $0 \leq ~ L_{f} <~\infty$. 
\end{myas}
\ras{lipschitzf} will be used to derive several conditions to guarantee feasibility of the OCP. In the formulation of MPC, we iteratively find at each update time $t_k$, $k\in \mathbb{N}_{\geq 0}$, an optimal predictive state $\hat{x}^* (\xi)$ and a control trajectory $u^* (\xi) $ for all $\xi \in [t_k , t_k + T_k]$, by minimizing the cost given by \req{cost}. 
Following the idea from \cite{Chen1998a}, we impose here the so called \textit{terminal constraint}, that the predictive state reaches $\Phi_f$ within $T_k$, i.e., $\hat{x}^* (t_k + T_k) \in \Phi_f$. 
Since $\hat{x}^* (t_k + T_k) \in \Phi_f$, there exists a positive time interval when the optimal state reaches the boundary of $\Phi_f$; let $T^* _k \in \mathbb{R}_{>0}$ be such time interval given by $\hat{x}^* (t_k +T^* _k )\in \partial {\Phi}_f$. The time instant $t_k + T^* _k$ is also illustrated in \rfig{terregion_robust}. 

Based on the above notations, we propose the following OCP:
\begin{myprob}[OCP] \label{OCP}
For the non-initial time $t_k$, $k\in \mathbb{N}_{\geq 1}$, 
given $x(t_k)$ and $T^* _{k-1}$, the OCP is to minimize the cost $J(x(t_k), u(\cdot ))$ given by \req{cost}, subject to 
\begin{numcases}
   { }
      \dot{\hat{x}}(\xi )=f (\hat{x}(\xi ) , u(\xi ) ),\ \xi \in [t_k, t_k+T_k] \label{constraint1} \\ 
    {u}(\xi ) \in {\mathcal U} \label{constraint2} \\
    \hat{x} (t_k + T_k )\in {\Phi}_{f}, \label{constraint3}
\end{numcases}
where $T_k = T^* _{k-1} - \alpha \Delta_{k-1}$ for a given $0< \alpha <1$ and $\Delta_{k-1} = t_k - t_{k-1}$. For the initial time $t_0$, minimize the cost $J(x(t_k), u(\cdot ))$ given by \req{cost}, subject to \req{constraint1}, \req{constraint2} and $\hat{x} (t_0 + T_0 )\in {\Phi}_{f}$ for a given $T_0 > 0$.
\end{myprob}

For the initial time $t_0$, Problem~1 is solved with a given $T_0 >0$. In order to guarantee the feasibility at $t_0$, $T_0$ needs to be suitably chosen such that the terminal constraint $\hat{x} (t_0 + T_0) \in \Phi_f$ is fulfilled. 
More specifically, $T_0$ should be selected to satisfy $x(t_0) \in {\mathcal X} (T_0)$, where ${\mathcal X} (T_0) = \{ x(t_0) \in \mathbb{R}^n\ |\ \exists u(t) \in {\mathcal U}, t\in [t_0, t_0+T_0]\ :\hat{x}(t_0 +T_0) \in {\Phi}_f\}$, i.e.,  ${\mathcal X} (T_0)$ denotes the set of states that can reach ${\Phi}_f$ within the time $t_0+T_0$. 
Although there may not exist a general framework to compute ${\mathcal X} (T_0)$ explicitly for nonlinear systems, several approximation methods have been proposed to compute ${\mathcal X} (T_0)$, see e.g., \cite{reachability2008}. 
The initial feasibility is essentially required for guaranteeing recursive feasibility, which is analyzed in the next section. 

For the \textit{non}-initial time $t_k, k\in \mathbb{N}_{\geq 1}$, 
we require by \req{constraint3} that the optimal state enters $\Phi_f$ within $T_k = T^* _{k-1} - \alpha \Delta_{k-1}$, where $T^* _{k-1}$ is the time interval obtained by the previous calculation of OCP. This implies that $T^* _k$ satisfies $T^* _k \leq T_k = T^* _{k-1} - \alpha \Delta_{k-1} < T^* _{k-1} \leq T_{k-1}$, which guarantees that the time interval $T^* _k$ and the prediction horizon $T_k$ become strictly smaller than the previous one at $t_{k-1}$. In later sections, we will make use of this property to show that the state enters $\Phi$ in finite time. 
\begin{myrem}
\normalfont 
Although various analysis and control strategies have been proposed for MPC, approaches to guarantee stability can be mainly divided into two categories; the OCP with a terminal constraint (see e.g., \cite{Chen1998a}), and the OCP \textit{without} a terminal constraint (see e.g., \cite{andrea2014a,grune2009a}). While the OCP becomes in general harder to be solved when the terminal constraint is imposed, this paper follows the former approach to guarantee stability and to derive an event-triggered strategy. 
Note that our problem formulation slightly differs from the standard formulation \cite{Chen1998a}, since the prediction horizon is not constant but is adaptively selected for each calculation time of the OCP. 
\qedwhite 
\end{myrem}

\section{Feasibility analysis} 
The main focus of this section is to derive several conditions to guarantee the notion of \textit{recursive feasibility}, which states that the existence of a solution to Problem~1 at an initial update time $t_0$ implies the feasibility at any update times afterwards $t_k$, $k\in \mathbb{N}_{\geq 1}$. 
The obtained feasibility conditions will be key ingredients to derive the event-triggered strategy, 
which will be discussed in the next section. 
\begin{mythm}\label{robust_feasibility}
Suppose that the OCP defined in \rprob{OCP} has a solution at $t_k$, providing an optimal control input $u^* (\xi)$ and the corresponding state trajectory $\hat{x}^* (\xi)$ for all $\xi \in [t_k, t_k + T_k]$, and the time $t_k + T^* _k$. 
Then, \rprob{OCP} has a solution at $t_{k+1} (>t_k)$, if the followings are satisfied: 
\begin{numcases}
   { }
||x (t_{k+1} ) - \hat{x}^* (t_{k+1} ) ||_{P_f} \leq (\varepsilon - \varepsilon_f ) e^{-L_{f}  T^* _k} \label{feasibility_cond} \\ 
\Delta_k =t_{k+1} - t_k \leq T^* _k \label{Deltak}, \\
   ||w (t) ||_{P_f} \leq \tilde{w} _{\rm max},\ \ \forall t \in [t_k, t_{k+1}] \label{disturbance_cond}
\end{numcases}
where $\tilde{w} _{\rm max}= \frac{\lambda_{\rm min} (\hat{Q}_{P}) }{ 4 e^{L_{f}  T^* _0} } (1-\alpha) \varepsilon_f$.
\end{mythm}
\begin{pf}
Consider the following dual mode controller as a feasible control candidate: 
\begin{equation}\label{feasibleinput_robust}
\begin{array}{lll}
\bar{u} (\xi )  = \left \{
\begin{array}{l}
u^* (\xi ),\ \  \xi \in[t_{k+1}, t_k +T^* _k ] \\
\kappa (\bar{x}(\xi )), \ \ \xi \in ( t_k +T^* _k, t_{k+1} + T_{k+1} ],
\end{array}
\right.
\end{array}
\end{equation}
where $T_{k+1} = T^* _k - \alpha \Delta_k$ and $\bar{x} (\xi)$ denotes the predictive state trajectory obtained by applying $\bar{u} (\xi )$, i.e., $\dot{\bar{x}} (\xi ) = f (\bar{x} (\xi ), \bar{u}(\xi ))$ with $\bar{x}(t_{k+1}) = x(t_{k+1})$. Note that we have $t_{k+1} + T_{k+1} > t_k +T^* _k$ since $t_{k+1} + T_{k+1} = t_k + \Delta_k + T^* _k - \alpha \Delta_k
 = t_k + (1-\alpha) \Delta_k + T^* _k > t_k + T^* _k$.
Furthermore, we have $T_{k+1} >0$ since $T^* _k -\alpha \Delta_k \geq (1-\alpha ) T^* _k >0$ from the condition \req{Deltak}. 

To prove that \req{feasibleinput_robust} is a feasible controller for $t_{k+1}$, we show that the following two arguments are satisfied:
\begin{enumerate}
\item By applying $\bar{u} (\xi )$, $\xi \in [t_{k+1}, t_k + T^* _k]$, the predictive state enters $\Phi$ by the time $t_k + T^* _k$. That is, $\bar{x} (t_k +T^* _k) \in \Phi$. This ensures that applying the local controller $\kappa$ from $t_k + T^* _k$ is admissible.  
\item By applying $\bar{u} (\xi )$, $\xi \in (t_{k}+T^* _k, t_{k+1} + T_{k+1}]$, the predictive state $\bar{x}$ enters ${\Phi}_f$ by the time $t_{k+1} + T_{k+1}$. That is, $\bar{x} (t_{k+1} + T_{k+1} )  \in {\Phi}_f$.
\end{enumerate}
The basic idea is to derive the upper bound of difference between $\bar{x}$ and $\hat{x}^*$ and show that the difference is small enough to prove (i), (ii); the reader can also refer to \cite{mpc_recent1} for a related analysis. To prove (i), we first use the Gronwall-Bellman inequality \cite{khalil} to obtain the upper bound of the difference between $\bar{x}$ and $\hat{x}^*$; 
$|| \bar{x} (\xi ) - \hat{x}^* (\xi ) ||_{P_f} \leq  || x(t_{k+1}) - \hat{x}^* (t_{k+1})||_{P_f} e^{L_{f}  (\xi -t_{k+1})}$ 
 for $\xi \in [t_{k+1}, t_{k} +T^* _k]$. 
Supposing that \req{feasibility_cond} holds and by letting  $\xi = t_k + T^* _k$, we obtain $|| \bar{x} ( t_k + T^* _k ) - \hat{x}^* ( t_k + T^* _k ) ||_{P_f} \leq e^{-L_{f}  T^* _k} (\varepsilon - \varepsilon_f ) e^{L_{f}  (t_k + T^* _k -t_{k+1})} = (\varepsilon - \varepsilon_f ) e^{-L_{f}  (t_{k+1}-t_k)}$.
Thus, from the triangle inequality, we obtain 
\begin{equation*}
\begin{array}{lll}
|| \bar{x} (t_k + T^* _k)||_{P_f} \\
\ \ \ \ \leq || \hat{x}^* (t_k + T^* _k) ||_{P_f} + (\varepsilon - \varepsilon_f ) e^{-L_{f}  (t_{k+1}-t_k)} \\
\ \ \ \ \leq  \varepsilon_f + \varepsilon - \varepsilon_f =  \varepsilon.
\end{array}
\end{equation*}
Thus it holds that $\bar{x} (t_k + T^* _k)\in \Phi$ and the proof of (i) is completed. 


We now prove the statement in (ii). 
By using $\bar{x} (t_k + T^* _k)\in \Phi$ and from \rlem{lem1}, we obtain $\dot{V}_f (\bar{x} (\xi )) \leq  -0.5 \ \bar{x} ^\mathsf{T} (\xi ) (Q+K^\mathsf{T} R K) \bar{x} (\xi ) \leq - 0.5\ \lambda_{\rm min} (\hat{Q}_{P}) V_f ( \bar{x} (\xi ))$ 
for $\xi \in (t_k + T^* _k, t_{k+1}+T^* _k - \alpha \Delta_k]$. 
Furthermore, by again applying the Gronwall-Bellman inequality, we obtain
\begin{equation*}
\begin{array}{lll}
||\bar{x} (t_k +T^* _k) ||_{P_f} \\
\ \ \ \leq  ||\hat{x}^* (t_k +T^* _k)||_{P_f} + \cfrac{\tilde{w} _{\rm max}}{L_{f} } e^{L_{f}  T^* _k} (1-e^{-L_{f}  \Delta_k}) \\
\ \ \ \leq  \varepsilon _f + \cfrac{(1-\alpha)}{4 L_{f} } \varepsilon_f  \lambda_{\rm min} (\hat{Q}_{P})(1-e^{-L_{f}  \Delta_k}),
\end{array}
\end{equation*}
where we have used \req{disturbance_cond} in the last inequality. Denoting $\eta = \frac{(1-\alpha)}{4 L_{f} } \lambda_{\rm min} (\hat{Q}_{P})$, and by using comparison lemma \cite{khalil}, we obtain 
$V_f ( \bar{x} (t_{k+1} + T^* _k -\alpha \Delta_k)) \leq   \varepsilon^{2} _f \left ( 1 + \eta (1-e^{-L_{f}  \Delta_k}) \right )^2  e^{-2 L_{f}  \eta \Delta_k }\leq  \varepsilon^2 _f$. 
The second inequality follows from the fact that the function $g_\varepsilon (\Delta_k) = ( 1 + \eta (1-e^{-L_{f}  \Delta_k}) )  e^{-L_{f}  \eta  \Delta_k } $ is shown to be a decreasing function of $\Delta_k$ with $g_\varepsilon (0) = 1$. 
Thus we obtain $V_f (\bar{x} (t_{k+1} + T^* _k -\alpha \Delta_k)) \leq \varepsilon^2 _f$, and the proof of (ii) is completed. Based on above, \req{feasibleinput_robust} is proven to be a feasible controller for $t_{k+1} (>t_k)$, provided that \req{feasibility_cond}, \req{Deltak} and  \req{disturbance_cond} are satisfied. This completes the proof. \qedwhite  
\end{pf}

\section{Event-triggered strategy}
By making use of the feasibility conditions provided in the previous section, we now propose an event-triggered strategy. Suppose again that the OCP is solved at $t_k$, providing a pair of optimal control input $u^* (\xi)$ and the corresponding state trajectory $\hat{x}^* (\xi)$ for all $\xi \in [t_k, t_k + T_k]$. 
In the following, event-triggered conditions based on the feasibility result will be derived to determine the next calculation time of the OCP $t_{k+1} (> t_k)$. 

The simplest way to determine $t_{k+1}$ might be to use the original feasibility conditions directly as the event-triggered conditions, i.e., for each $t > t_k$, check the feasibility according to \req{feasibility_cond} and \req{Deltak}, i.e., 
\begin{equation}\label{ev_condition}
||x (t) - \hat{x}^* (t ) ||_{P_f} \leq (\varepsilon - \varepsilon_f ) e^{-L_{f}  T^* _k},
\end{equation}
\begin{equation}\label{ev_condition2}
t - t_k \leq T^* _k.
\end{equation}
Only when either of the above conditions is violated, then we set $t_{k+1} = t$ as the next update time. 
However, checking the above conditions for each $t > t_k$ requires \textit{continuous} monitoring of the state $x(t)$ and evaluation of the above conditions, 
which may lead to a high cost of sensing requiring a dedicated analog hardware, and thus it is not suitable for standard digital platforms used in real-time implementations. 

\begin{figure}[tbp]
  \begin{center}
   \includegraphics[width=6.0cm]{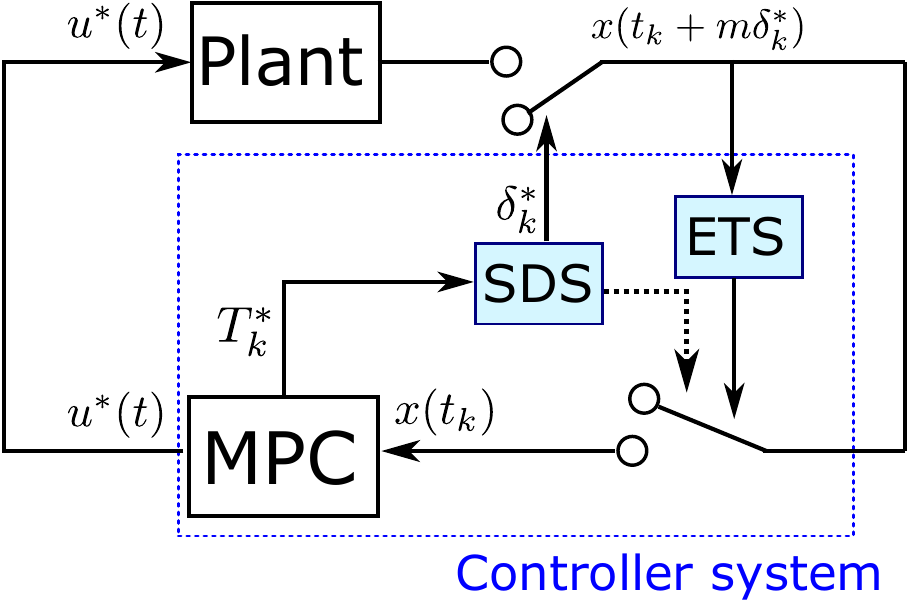}
   \caption{ A schematic overview of the event-triggered strategy.} 
   \label{system_configuration}
  \end{center}
 \end{figure}

Therefore, we propose here an alternative event-triggered approach by relaxing the above \textit{continuous} requirements.
The key idea of our approach is to measure the state and evaluate event-triggered conditions only at certain sampling time instants, instead of continuously. A schematic block diagram of our proposed scheme is illustrated in \rfig{system_configuration} and the overview is stated as follows. 
Once the OCP is solved by MPC at an update time instant, say $t_k$, and $T^* _k$ is obtained, 
the Sampling-time Decision System (SDS) computes $\delta^* _k \in \mathbb{R}_{>0}$, which represents the sampling time interval at which the event-triggered condition is evaluated. Namely, from the obtained $\delta^* _k$ from SDS, the Event-Triggered System (ETS) measures the state and checks the event-triggered condition only at $t_k + m\delta^* _k$, $m\in \mathbb{N}_{\geq 1}$, in order to determine the next update time $t_{k+1}$. 
Note that the SDS has a partial role to determine $t_{k+1}$ to solve the OCP (the black dotted arrow in \rfig{system_configuration}); as described later in this section, $t_{k+1}$ can sometimes be directly determined according to $T^* _k$ without needing to evaluate the event-triggered condition. 

\begin{figure}[tbp]
  \begin{center}
   \includegraphics[width=6.0cm]{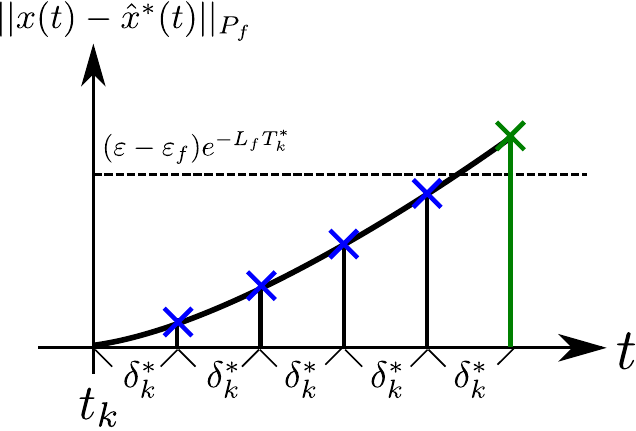}
   \caption{The figure illustrates the problem of violating the feasibility described in (P.2). The blue marks represent the sequence of the left hand side in \req{ev_condition} evaluated with the sampling time $\delta^* _k$. 
As shown in the figure, the feasibility can be violated between two evaluation times (represented as green mark).}
   \label{periodic_illustration}
  \end{center}
 \end{figure}

Regarding the proposed framework outlined above, we need to derive both mechanisms to determine $\delta^* _k$ and the event-triggered conditions. One might directly utilize \req{ev_condition}, \req{ev_condition2} as the event-triggered conditions, and evaluate them with a given arbitrary value of $\delta^* _k$. 
However, this cannot be applied due to the following two problems regarding the violation of feasibility:


\renewcommand{\labelenumi}{(\textbf{P.\arabic{enumi})}}
\begin{enumerate}
\item If a large value of $\delta^* _k$ would be chosen, the feasibility would not be satisfied (i.e., the left hand side of \req{ev_condition} exceeds the threshold in the right hand side)  at the next evaluation time $t_k + \delta^* _k$. 
\item If we would directly use \req{ev_condition} as the event-triggered condition, the feasibility might be violated between two consecutive evaluation times. This issue is illustrated in \rfig{periodic_illustration}. 
\end{enumerate}

In the following, we provide solutions to each problem above and then provide the overall event-triggered strategy. 
Consider first to solve (P.1); to deal with the problem, $\delta^* _k$ needs to be chosen small enough such that the feasibility conditions \req{ev_condition}, \req{ev_condition2} are satisfied for all $t \in [t_k, t_k + \delta^* _k]$. Thus, let us consider to evaluate a minimum inter-event time of the feasibility conditions \req{ev_condition}, \req{ev_condition2}. Assume that the size of the disturbance satisfies $||w(t)||_{P_f} \leq \tilde{w} _{\rm max}, \forall t \geq t_0$, which ensures from \rthm{robust_feasibility} that the effect of disturbances does not violate the feasibility. 
By using Gronwall-Bellman inequality, we obtain 
$||x(t) - \hat{x}^* (t)||_{P_f} \leq \frac{\tilde{w} _{\rm max}}{L_{f} } (e^{L_{f}  (t-t_k)}-1)$ 
for $t\in [t_k, t_k + T_k]$. Thus, a sufficient condition to satisfy \req{ev_condition} is 
$\frac{\lambda_{\rm min} (\hat{Q}_{P} ) (1-\alpha) \varepsilon_f}{ 4L_{f}  e^{L_{f}  T^* _0}}  ( e^{L_{f}  (t-t_k)}-1) \leq (\varepsilon - \varepsilon_f ) e^{-L_{f}  T^* _k}$.
Solving the inequality for $t$ yields 
$t \leq t_k + \Delta^{\rm min} _k$, where 
\begin{align}\label{delta_minimum}
\Delta^{\rm min} _k =\cfrac{1}{L_{f} } \ln \left (1 + \cfrac{ 4 L_{f} (\varepsilon - \varepsilon_f ) e^{L_{f}  (T^* _0 -T^* _k )}}{\lambda_{\rm min} (\hat{Q}_{P} ) (1-\alpha) \varepsilon_f }  \right )
\end{align}
Thus, the condition \req{ev_condition} is satisfied for all $t\in [t_k, t_k + \Delta^{\rm min} _k]$. By taking into account the other feasibility condition \req{ev_condition2}, the over-all minimum inter-event time is now given by $\min \{ \Delta^{\rm min} _k,\ T^* _k \}$. 
For the case when $\Delta^{\rm min} _k \leq T^* _k$, the minimum inter-event time becomes $\Delta^{\rm min} _k$. 
Thus, if $\delta^* _k$ is selected such that $\delta^* _k = \gamma \Delta^{\rm min} _k \leq \Delta^{\rm min} _k$ for a given $0 < \gamma \leq 1$, the feasibility condition is fulfilled for all $t \in [t_k , t_k + \delta^* _k ]$. 
On the other hand, for the case when $\Delta^{\rm min} _k > T^* _k$, \req{ev_condition} is satisfied for all $t \in [t_k , t_k + T^* _k]$. This means that \req{ev_condition2} is violated earlier than \req{ev_condition}. 
Thus, if $T^* _k < \Delta^{\rm min} _k$, the next update time can be directly set as $t_{k+1} = t_k + T^* _k$. 
Based on the above analysis, the following strategy can be provided as a solution to (P.1): 
\renewcommand{\labelenumi}{(\alph{enumi})}
\begin{enumerate}
\item If $T^* _k \geq \Delta^{\rm min} _k$, then set $\delta^* _k = \gamma \Delta^{\rm min} _k$ for a given $0< \gamma \leq 1$. 
\item If $T^* _k < \Delta^{\rm min} _k$, then set $t_{k+1} = t_k + T^* _k$ as the next update time. 
\end{enumerate}

Note that the above strategy is implemented by the SDS, i.e., depending on the value of $T^* _k$, it either determines $\delta^* _k$ (a) or directly the next update time $t_{k+1}$ (b). 
From \req{delta_minimum}, $\Delta^{\min} _k$ as well as $\delta_k^* ( = \gamma \Delta^{\rm min} _k)$ get larger as $T^* _k$ decreases. Therefore, if the initial time fulfills $T^* _0 > \Delta^{\min} _0$, case (a) is selected for the initial time and both $\Delta^{\min} _k$ and $\delta^* _k$ increase afterwards while $T^* _k$ decreases. 
Since $\Delta^{\min} _k$ gets larger and $T^* _k$ decreases, the magnitude relation will be switched to $T^* _k < \Delta^{\min} _k$ (case (b)) after a certain time step. As shown in case (b), the computation of $\delta^* _k$ is no more required for this case since the next update time $t_{k+1}$ can be directly determined as $t_{k+1} = t_k + T^* _k$. 
\begin{myrem}
\normalfont 
One may argue that the inter-event time is given by $T^* _k$ for case (b) and that it may thus tend to $0$ since $T^* _k$ is decreasing. Note however, that $T^* _k > 0$ always holds while the MPC is implemented (i.e., $x(t_k) \notin \Phi$); if $x(t_k)$ is outside of $\Phi$, there always exists a strictly positive time interval for the optimal state to reach $\Phi_f$. 
Thus, this guarantees that the inter-event time remains always positive while implementing the MPC. 
\qedwhite 
\end{myrem}


Next, consider the problem (P.2). 
Based on the obtained $\delta^* _k$, the time instants to measure the state and evaluate the event-triggered condition are now given by $t_k + m \delta^* _k$, $m \in \mathbb{N}_{\geq 1}$. To avoid losing the feasibility between two evaluation times, the ETS checks the feasibility condition at \textit{one time step ahead}, instead of the current time instant. That is, at an evaluation time $t = t_k + m\delta^* _k$, $m \in \mathbb{N}_{\geq 1}$, 
the feasibility is checked for $t + \delta^* _k$ instead of $t$. If the feasibility at $t+ \delta^* _k$ \textit{will be} guaranteed, then the ETS moves on to the next evaluation time $t+ \delta^* _k$, and repeats the same procedure. 
On the other hand, if the feasibility at $t + \delta^* _k$ is not guaranteed, then the next update time is set as $t_{k+1} = t$. 
Since the ETS preliminary checks the feasibility at one step future time, the loss of feasibility does not occur between two evaluation times. 

The feasibility at one time step ahead can be given by modifying the original feasibility conditions. Suppose at an evaluation time $t = t_k + m\delta^* _k$, $m \in \mathbb{N}_{\geq 1}$, we aim at checking the feasibility at $t+\delta^* _k$ based on the current state $x(t)$. 
The difference between the actual state and the optimal state at $t + \delta^* _k$ is given by
\begin{equation}\label{next_state_difference}
\begin{array}{lll}
||x (t + \delta^* _k ) - \hat{x}^* (t + \delta^* _k) ||_{P_f} \\ 
\ \ \ \ \ \ \leq  e^{L_{f}  \delta^* _k} ||x (t) - \hat{x}^* (t) ||_{P_f} + \cfrac{\tilde{w}_{\max} }{ L_f }(e^{L_f \delta^* _k} -1 ) 
\end{array}
\end{equation}
where we have used \req{disturbance_cond}. 
From the feasibility conditions \req{ev_condition}, \req{ev_condition2}, feasibility at $t + \delta^* _k$ is guaranteed if $t+ \delta^* _k - t_k \leq T^* _k$ and $||x (t+ \delta^* _k) - \hat{x}^* (t+ \delta^* _k) ||_{P_f}< (\varepsilon - \varepsilon_f ) e^{-L_{f}  T^* _k}$ are both satisfied. 
From \req{next_state_difference}, sufficient conditions to satisfy these conditions are thus given by
\begin{align}
||x & (t) - \hat{x}^* (t) ||_{P_f} \notag \\
           &< (\varepsilon - \varepsilon_f ) e^{-L_{f}  (T^* _k + \delta^* _k)} - \cfrac{\tilde{w}_{\max} }{ L_f }(1-e^{-L_f \delta^* _k}) \label{periodic_ev_condition} \\ 
 & (m+1) \delta^* _k \leq T^* _k. \label{periodic_ev_condition2}
\end{align}
Note that if \req{periodic_ev_condition},  \req{periodic_ev_condition2} are both satisfied the feasibility is guaranteed at $t + \delta^* _k$, 
and these conditions can be evaluated based on $x(t)$. 
Therefore, by letting the ETS evaluate \req{periodic_ev_condition} and \req{periodic_ev_condition2} as the event-triggered conditions, the violation of the feasibility between two evaluation times will not occur, providing thus a solution to (P.2). 

The over-all proposed algorithm of the event-triggered strategy is now summarized below: 

\noindent
{\bf {Algorithm 1} (Event-triggered MPC)}:
\renewcommand{\labelenumi}{(\roman{enumi})}
\begin{enumerate}
\item At any update times $t_k$, $k\in \mathbb{N}_{\geq 0}$, if $x(t_{k}) \in \Phi$, then switch to the local controller $\kappa (x)$ as a dual mode strategy. Otherwise, solve \rprob{OCP} and obtain the optimal control and state trajectory $u^* (\xi )$,  $\hat{x}^* (\xi )$ for all $\xi \in[t_k, t_k +T_k]$, and $T^* _k$ as the time interval when the state reaches $\Phi_f$, i.e., $\hat{x}^* (t_k + T^* _k) \in \partial \Phi_f$. 
\item The SDS provides the sampling time $\delta^* _k$ or the next update time $t_{k+1}$ in the following way: 
\begin{enumerate}
\item If $T^* _k \geq \Delta^{\rm min} _k$, 
then set $\delta^* _k = \gamma \Delta^{\rm min} _k$ for a given $0 < \gamma \leq 1$, and go to step (iii). 
\item If $T^* _k < \Delta^{\rm min} _k$, then set 
$t_{k+1} = t_k + T^* _k$ and go to step (iv). 
\end{enumerate}
\item The ETS provides the next update time $t_{k+1} (>t_k)$ in the following way:
\begin{enumerate}
\item Set $m=1$. 
\item At an evaluation time $t = t_k + m \delta^* _k$, $m\in \mathbb{N}_{\geq 1}$, measure the state $x(t)$, and check the event-triggered conditions given by \req{periodic_ev_condition}, \req{periodic_ev_condition2}.
\item If \req{periodic_ev_condition} and \req{periodic_ev_condition2} are both satisfied, then apply $u^* (\xi)$ for $\xi \in [t, t + \delta^* _k )$. Then, set $m \leftarrow m+1$ and go back to step (b). 
Otherwise, set $t_{k+1} = t $ and go to step (iv). 
\end{enumerate}
\item $k \leftarrow k+1$ and go back to step (i). \qedwhite
\end{enumerate}

\section{Stability analysis}
In this section we analyze stability of the closed loop system under the implementation of Algorithm~1.  
We will prove in the following that, any state trajectories starting from the initial feasible set ${\mathcal X} (T_0)$ (see the definition of ${\cal X} (T_0)$ in Section~2) will eventually enter $\Phi$ within a prescribed finite time interval. 

\begin{mythm}\label{convergence}
Consider the nonlinear system given by \req{sys1}, and suppose that Algorithm~1 is implemented. 
Then, for any $w(t)$ satisfying $||w(t)||_{P_f} \leq \min \{\hat{w} _{\max}, \tilde{w} _{\rm max} \}$, $\forall t \geq t_0$, any state trajectories starting from $x(t_0) \in {\mathcal X}(T_0)$ enter $\Phi$ within the time interval $T^* _0 /\alpha$, and remain in $\Phi$ for all the future times. 
\end{mythm}
\begin{pf}
We prove the statement by contradiction. Assume that at $t_k$ we have $t_k - t_0 \geq T^* _0 /\alpha$, and $x(t_k)$ is outside of $\Phi$, i.e., $x(t_k) \notin \Phi$. 
Since $x(t_k) \notin \Phi$ and $\Phi_f \subset \Phi$, we have $T^* _k > 0$. 
As $x(t_0) \in {\mathcal X}(T_0)$ and $||w(t)||_{P_f} \leq \tilde{w} _{\rm max}$, $\forall t \geq t_0$, applying Algorithm 1 ensures that the feasibility is guaranteed for all $t_0, t_1, \cdots, t_k$. Thus, we recursively obtain from \req{constraint3} that: 
\begin{equation*}
\begin{aligned}
T^* _k &\leq T^* _{k-1} - \alpha \Delta_{k-1} \leq T^* _{k-2} - \alpha ( \Delta_{k-1} + \Delta_{k-2}) \\
         &\leq \cdots \leq T^* _0 - \alpha \sum^{k-1} _{l=1} \Delta_l \\
& = T^* _0 - \alpha (t_k -t_{k-1} + t_{k-1} - t_{k-2} + \cdots + t_1 -t_0) \\ 
&= T^* _0 - \alpha (t_k - t_0).
\end{aligned}
\end{equation*}
Thus, by the assumption $t_k - t_0 \geq T^* _0 /\alpha$, we obtain $T^* _k  \leq 0 $. However, this clearly contradicts the fact that we have $T^* _k > 0$. Thus, it is shown that the state enters $\Phi$ within the time interval $T^* _0/\alpha$. 
Furthermore, since from \rlem{lem1}, $\Phi$ is a positively invariant set with the disturbance satisfying $||w(t)||\leq \hat{w} _{\rm max}$, 
the state remains in $\Phi$ for all future times. 
This completes the proof. \qedwhite
\end{pf}

\begin{myrem}[\textbf{On the control performance\\}]
\normalfont 
In~\rthm{convergence}, stability is proven by evaluating a time interval to reach $\Phi_f$, and not by the optimal cost. 
Although this may be unconventional with respect to a control performance view point, our approach is advantageous and practical from a \textit{event-triggered control} view point, since, as previously mentioned in Section~1, the event-triggered condition provides less conservative results than the existing ETMPC approaches. Moreover, the control performance can be evaluated by tuning the parameter $\alpha$. For more details, please see \rrem{rem1} below. \qedwhite 
\end{myrem}

\begin{myrem}[Convergence time v.s. Disturbance]\label{rem1}
\normalfont
If $\alpha$ is chosen larger, then $T^* _0 / \alpha$ gets smaller and faster convergence is obtained. However, this in turn means from \req{disturbance_cond} that the allowable size of disturbance becomes smaller, which implies that the robustness to noise or model uncertainty may be degraded. Thus, there exists a trade-off between the convergence time of the state trajectory and the allowable size of the disturbance, and this trade-off can be regulated by tuning $\alpha$. \qedwhite
\end{myrem}

\section{Simulation Results}
\begin{figure}[tbp]
  \centering
   \includegraphics[width=7.5cm]{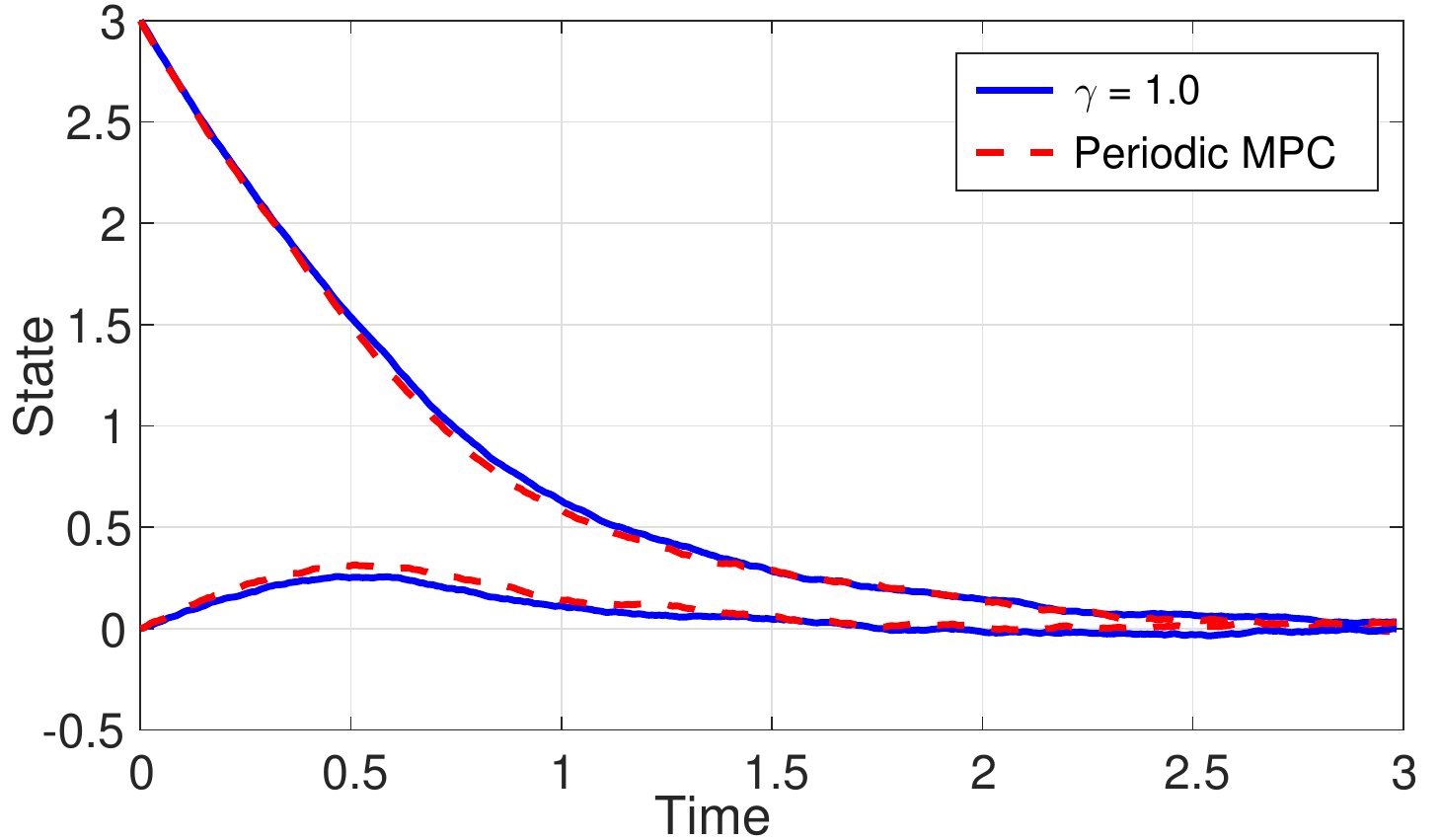}
   \caption{State trajectories under Algorithm~1 ($\gamma =1$) and the periodic MPC with a constant sampling time interval $0.1$.} 
   \label{state_traj_periodic}
 \end{figure}

As a simulation example, we consider the following system adopted from \cite{Chen1998a}: 
\begin{equation}
\begin{aligned}
\dot{x}_1 &= x_2 + u (\mu + (1-\mu ) x_1) + w_1 \\
\dot{x}_2 &= x_1 + u (\mu - 4 (1-\mu ) x_2) + w_2,
\end{aligned}
\end{equation}
with $\mu = 0.8$, where $x = [x_1; x_2] \in \mathbb{R}^2$, $u\in \mathbb{R}$, and $w = [w_1; w_2]\in \mathbb{R}^2$. We assume ${\cal U} = \{ u\in \mathbb{R}| -2 \leq u \leq 2 \}$, 
the matrices for the stage cost are $Q=0.1 I_2$, $R=0.05$, and the initial prediction horizon is set to $T_0 = 4.0$. The local controller is given by $\kappa = Kx$ with $K = [1.8042 \ 1.8042]$, and $P_f = [0.0814\ \ 0.0314;\ \ 0.0314 \ \ 0.0814]$ by following the procedure presented in \cite{Chen1998a}.
The computed Lipschitz constant is $L_{f} = 0.53$ and we set $\varepsilon_f = 0.08$, $\alpha = 0.8$. From \rthm{robust_feasibility}, the feasibility is guaranteed if $\tilde{w} _{\rm max} = 8.3 \times 10^{-4}$ and from \rlem{lem1} the region $\Phi$ is positively invariant if $\hat{w} _{\rm max} =2.0 \times 10^{-3}$. Taking into account both restrictions, we assume that ${\cal W} = \{ w \in \mathbb{R}^2 \ |\ ||w||_{P_f} \leq 8.3 \times 10^{-4} \}$. 
\rfig{state_traj_periodic} represents state trajectories of $x$ under Algorithm~1 with $x(t_0) = x(0) = [3; 0]$, $\gamma=1.0$, and the standard periodic MPC with a constant sampling time interval $0.1$. 
From the figure, the state trajectory under Algorithm~1 converges to a region around the origin similarly to the periodic case. 
The resulting convergence time needed for the state to enter $\Phi$ under the proposed method is $2.89 (\leq T^* _0 /\alpha = 3.88)$, and thus \rthm{convergence} is verified. 

\begin{figure}[tbp]
  \centering
   \includegraphics[width=7.5cm]{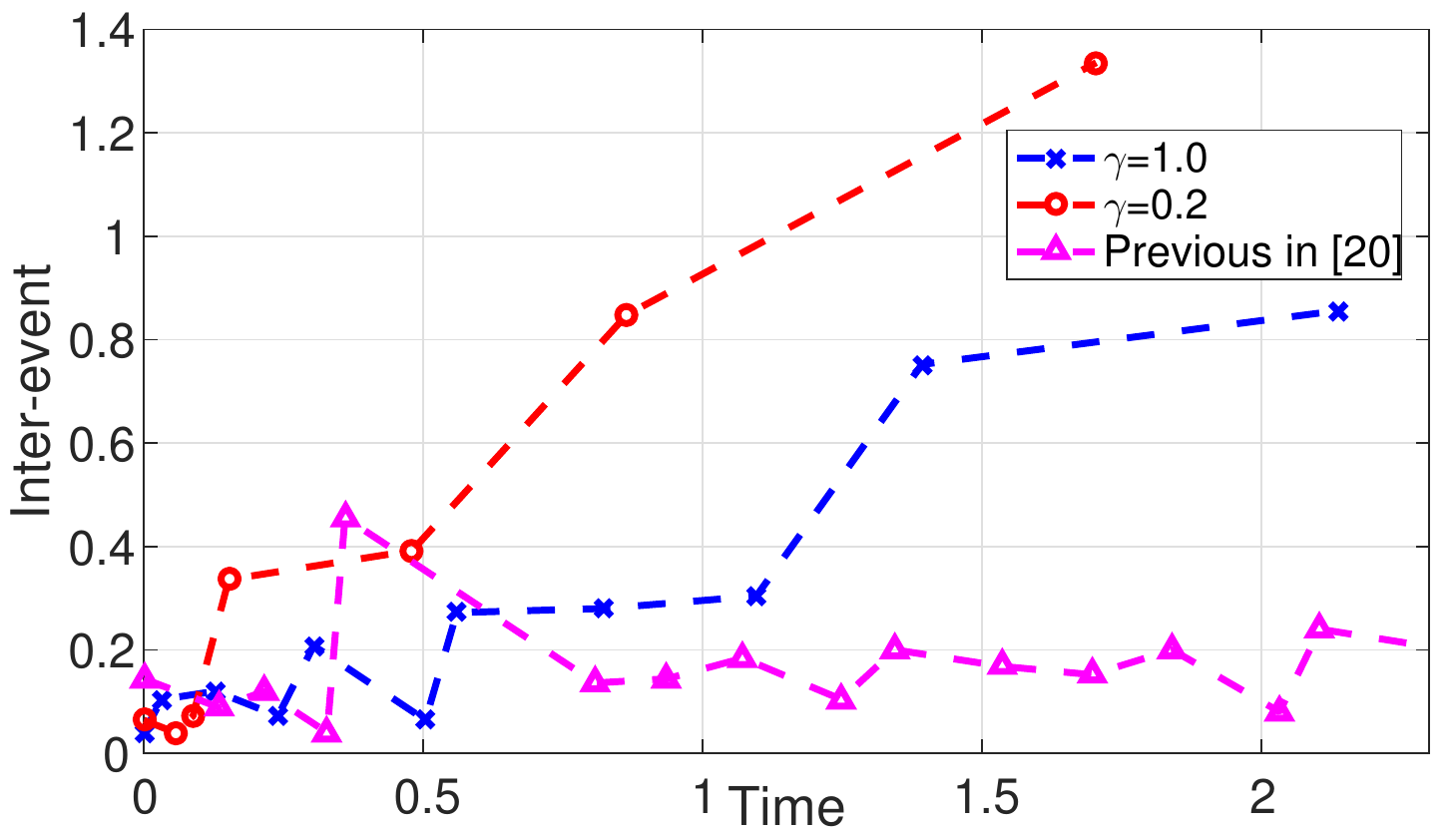}
   \caption{Inter-event times of solving OCPs under Algorithm~1 ($\gamma=0.2,\ 1$) and the previous approach presented in \cite{hashimoto2017a}.} 
   \label{triggering_events2}
 \end{figure}

\rfig{triggering_events2} shows the inter-event times $\Delta_k$ under Algorithm~1 ($\gamma = 0.2, 1.0$), and the approach presented in \cite{hashimoto2017a}. 
From the figure, the proposed scheme is shown to be more practical than our previous approach, since it achieves longer inter-event times. 

From the result, we can also evaluate the sensing cost by counting how often states are measured to check the event-triggered conditions. The total number of time instants when states are measured are given by $85$ (times) when $\gamma = 0.2$ and $46$ (times) when $\gamma=1.0$. Thus, less sensing cost is attained when $\gamma=1.0$. 
On the other hand, since the inter-event times are longer for the case when $\gamma=0.2$ than when $\gamma=1.0$ as illustrated in \rfig{triggering_events2}, less computational cost of solving the OCPs is achieved in the former case. 
Therefore, it is shown that there exists a trade-off between the sensing cost and the computational cost, and the trade-off can be regulated by tuning the parameter $\gamma$. 
 
In the simulation example, the allowable size of disturbance is given by $\tilde{w} _{\rm max} = 8.3 \times 10^{-4}$, which may be relatively small. Although this size may get larger by parameter tuning, how this restriction can be relaxed should be studied in our future research. 
\section{Conclusions}
In this paper, we proposed an event-triggered strategy for MPC of nonlinear continuous-time systems with additive bounded disturbances. 
The proposed method is derived based on new feasibility and stability results by imposing the terminal constraint with an adaptive prediction horizon. In the derivations of stability, we evaluate the time interval when the optimal state trajectory enters the local set $\Phi_f$, and it is shown that the state converges $\Phi$ within a prescribed finite time interval. Furthermore, the proposed event-triggered conditions are evaluated only at certain sampling time instants, aiming at reducing sensing cost and are thus suitable for practical implementations. A simulation example illustrates the effectiveness of the proposed scheme. 

\appendix
\section{ Proof of Lemma 1}
 Consider a linearization of \req{sys1} around the origin for the non-disturbance case; 
$\dot{x}(t) = A_f x(t) + B_f u(t)$, 
where $A_f = \partial f /\partial x (0, 0)$ and $B_f = \partial f/ \partial u (0, 0)$. Since the linearized system is stabilizable from \ras{as1}, we can find a state feedback controller $\kappa(x) = K x$ such that $A_c = A_f +B_f K$ is Hurwitz and the closed loop system $\dot{x} = A_c x$ is thus asymptoptically stable. 
Choose a matrix $P$ such that the following Lyapunov equation holds: $P A_c + A^\mathsf{T} _c  P = - (Q+ K^\mathsf{T} R K)$ where $Q$ and $R$ are matrices for the stage cost defined in \req{cost}. Then, the time derivative of the function $V_f = x^\mathsf{T} P x$ along a trajectory of the nominal system $\dot{x} = f(x, \kappa(x))$ yields:
\begin{equation*}
\begin{aligned}
\dot{V}_f & ( x) = - x^\mathsf{T} (Q+K^\mathsf{T} R K) x + 2x^\mathsf{T} P \phi (x) \\
                 \leq & - x^\mathsf{T} (Q+K^\mathsf{T} R K) x \left ( 1-  \cfrac{2||\phi(x)||_P}{\lambda_{\rm min} (\hat{Q}_P) ||x||_P} \right), 
\end{aligned}
\end{equation*}
where $\phi (x) = f(x, \kappa(x)) - A_c x$, and $\hat{Q}_{P} =P^{-1/2} ( Q+K^\mathsf{T} R K )P^{-1/2}$. Since $||\phi(x)||_P / ||x||_P \rightarrow 0$ as $||x||_P \rightarrow 0$, there exists a positive constant $0 < \varepsilon_0 < \infty$ such that $||\phi(x)||_P / ||x||_P \leq \lambda_{\rm min} (\hat{Q}_P))/4$ for $||x||_P \leq \varepsilon_0$.
Let $0< \varepsilon \leq  \varepsilon_0$ such that for all $||x||_P \leq \varepsilon$, $\kappa(x) = K x \in {\mathcal U}$. By letting $\Phi = \{ x\in \mathbb{R}^n\ |\ V_f (x) \leq \varepsilon^2 \}$, we obtain $\dot{V}_f (x) \leq - 0.5 x^\mathsf{T} (Q+K^\mathsf{T} RK) x$ for all $x\in \Phi$. 

Now consider the time derivative of the function $V_f$ along a trajectory of the nonlinear system with additive disturbances $\dot{x} = f(x, \kappa(x)) + w$:
\begin{align*}
\dot{V}_f (x) &= - x^\mathsf{T} (Q+K^\mathsf{T} RK) x + 2x^\mathsf{T} P \phi (x) + 2 x^\mathsf{T} P w \notag \\
                                    &\leq - x^\mathsf{T} (Q+K^\mathsf{T} RK) x \left ( 1-  \cfrac{2 ||\phi(x)||_P}{\lambda_{\min } (\hat{Q}_P) ||x||_P} \right. \notag \\ 
& \ \ \ \ \ \ \ \ \ \ \ \ \ \ \ \ \ \ \ \ \ \ \ \  \left. - \cfrac{2 ||w||_P}{\lambda_{\min } (\hat{Q}_P)||x||_P} \right), \label{vdot_dist}
\end{align*}
and consider also a compact set as a boundary of $\Phi$; $\partial \Phi = \{ x\in \mathbb{R}^n\ |\ V_f (x) = \varepsilon^2 \}$. 
From above, we obtain $\dot{V}_f \leq 0$ for $x\in \partial \Phi$, if $||w||_P \leq \varepsilon\lambda_{\min} (\hat{Q}_P ) /4$.
Thus, $\Phi$ is a positive invariant set for the closed loop system $\dot{x} = f(x, \kappa(x)) + w$ if the disturbance satisfies $||w||_P \leq \varepsilon\lambda_{\min} (\hat{Q}_P )/4$. 
This completes the proof of Lemma 1.

\end{document}